\documentclass[11pt]{amsart}
\usepackage{hyperref}
\newtheorem{thm}{Theorem}[section]

\newtheorem{lem}[thm]{Lemma}

\theoremstyle{definition}

\theoremstyle{remark}
\newtheorem{rem}[thm]{Remark}
\numberwithin{equation}{section}

\begin{document}

\title[On the spectrum of the Laplace-Beltrami operator]{An application of scattering theory
\\ to the spectrum \\ of the Laplace-Beltrami operator}
\author{Francesca Antoci}
\address{Dipartimento di Matematica, Politecnico di Torino}
\email{antoci@dimat.unipv.it}

\subjclass{58G25, 35P25}
\keywords{Scattering theory; differential forms}%

\begin{abstract}
Applying a theorem due to Belopol'ski and Birman, we show that the
Laplace-Beltrami operator on 1-forms on ${\bf R}^n$ endowed with an
asymptotically Euclidean metric has absolutely continuous spectrum equal
to $[0, +\infty)$.
\end{abstract}
\maketitle
\section{Introduction}
The  relationships between the geometric properties of a complete
noncompact Riemannian manifold and the spectrum of the
Laplace-Beltrami operator have been intensively investigated by
many authors.\par  Among them, H. Donnelly since the late
seventies studied the spectra of the Laplacian and of the
Laplace-Beltrami operators on particular manifolds, such as the
hyperbolic space (\cite{Donnelly4}), manifolds with negative
sectional curvature (\cite{Donnelly5}), asymptotically euclidean
manifolds (\cite{Donnelly2}). However, to our knowledge, the case
of the Laplace-Beltrami operator acting on $p$-forms on
asymptotically euclidean manifolds has been left aside up to now.
\par The purpose of this paper is to contribute to the
investigation of this case. We study the absolutely continuous
spectrum of the Laplace-Beltrami operator on ${\bf R}^n$ endowed
with an asymptotically euclidean metric, that is with a Riemannian
metric satisfying conditions (\ref{1}), (\ref{2}) and (\ref{3}).
The tool employed is classical scattering theory in the wave
operators approach. In this case, the problem is reduced to a
problem of scattering for vector-valued operators.
\par An inherent restriction of the proof, however, that makes it
difficult an extension to more general cases, is the use of the
Fourier transform, which has a crucial role in our considerations,
particularly in connection with Lemma \ref{trace class}. The lack
of this tool in the more general case of manifolds with an
asymptotically controlled Riemannian metric is a major obstacle to
the extension of the theorem.\par

\section{Preliminaries}

Let $({\bf R}^n,e)$ be the euclidean $n$-dimensional space, and $({\bf
R}^n,g)$
 be the same space, endowed with a complete Riemannian metric $g$. \par
We will denote by $\Lambda^1_c({\bf R}^n)$ the vector space of all smooth,
compactly supported $1$-forms on  ${\bf R}^n$, and by $L^2_1({\bf R}^n,e)$
the completion of $\Lambda^1_c({\bf R}^n)$ with respect to the norm
\begin{equation} \label{L2e} \| \omega \|^2_{L^2_1({\bf R}^n,e)}= \int_{{\bf R}^n}
<\omega,\omega>_e \, dx,\end{equation} where $dx$ denotes the Euclidean
volume element and $$<\omega(x),\omega(x)>_e=\sum_i
\omega_i^2(x)$$ is the fiber norm for $1$-forms induced by the
Euclidean metric.
\par $L^2_1({\bf R}^n,e)$ is the Hilbert space
 direct sum of $n$ copies of $L^2({\bf R}^n)$. The
Laplace-Beltrami operator $\Delta_e$ on $1$-forms $\omega \in
\Lambda^1_c({\bf R}^n)$ acts componentwise as: $$ (\Delta_e \,
\omega)_k = -\sum_j \frac{\partial^2 \omega_k}{\partial x_j^2}.$$
It is well-known that $\Delta_e$ is essentially selfadjoint on
$\Lambda^1_c({\bf R}^n)$, and its closure $H_0$ has purely
absolutely continuous spectrum equal to $$\sigma(H_0)=[0,
+\infty),$$ with constant multiplicity. \par We will denote by
${\bf h}_0[\omega]$ the quadratic form associated to $\Delta_e$ on
$\Lambda_c^1({\bf R}^n)$:
\begin{equation}\label{h0} {\bf h}_0[\omega] = \int_{{\bf R}^n} <\Delta_e
\omega, \omega>_e\,dx= \int_{{\bf R}^n}\sum^n_{i,j=1} \left(\frac{\partial
\omega_i}{\partial x^j}\right)^2\, dx.\end{equation}

$L^2_1({\bf R}^n,g)$ will stand for the completion of
$\Lambda^1_c({\bf R}^n)$ with respect to the norm:
\begin{equation} \label{L2G} \|\omega\|^2_{L^2_1({\bf R}^n,g)}= \int_{{\bf R}^n}
<\omega,\omega>_g\,\sqrt{g}\,dx,\end{equation}
where $\sqrt{g}\, dx$, as usual, denotes the volume element
induced by the Riemannian metric $g$ and
$$<\omega(x),\omega(x)>_g=g^{ij}(x) \omega_i(x)\omega_j(x)$$ is
the fiber norm for $1$-forms induced by $g$. (Here, as everywhere
throughout the paper, the repeated indices convention is
adopted.)\par The action of the Laplace-Beltrami operator
$\Delta_g= d\delta + \delta d$  on $1$-forms $\omega \in
\Lambda^1_c({\bf R}^n)$  is given in local coordinates by the
Weitzenb\"ock formula: $$(\Delta_g \omega)_k= - (g^{ij}\nabla_i
\nabla_j \omega)_k + R^i_k \omega_i,$$ where $\nabla_i$ is the
covariant derivative with respect to the connection induced by the
the metric $g$, and $R^i_k$ is the Ricci tensor.
\par
Since the Riemannian metric $g$ is complete, $\Delta_g$ is
essentially selfadjoint on $\Lambda^1_c({\bf R}^n)$. We will
denote its closure by $H_1$.\par Moreover, we will denote by ${\bf
h}_1[\omega]$ the quadratic form associated to $\Delta_g$ on
$\Lambda_c^1({\bf R}^n)$
\begin{equation}\label{h1}
{\bf h}_1[\omega] = \int_{{\bf R}^n}<\Delta_g \omega, \omega>_g
\sqrt{g}\, dx= \end{equation} $$=\int_{{\bf R}^n} |\nabla \omega
|^2_g \sqrt{g}\, dx + \int_{{\bf R}^n}
<R \omega, \omega>_g \sqrt{g}\, dx, $$
where  $$|\nabla \omega |^2_g=g^{ij}g^{\alpha \beta} \nabla_i
\omega_{\alpha} \nabla_j \omega_{\beta} $$ and $$<R \omega,
\omega>_g=g^{\alpha \beta} R^i_{\alpha} \omega_i \omega_{\beta}.$$
\par In the next section we will show how it is possible, under
suitable hypothesis on the asymptotic behaviour of $g$, to get
information about the spectrum of $H_1$ from the knowledge of the
spectrum of $H_0$, proving the following
\begin{thm}\label{euclidean}
Let ${\bf R}^n$ be endowed with a Riemannian metric $g$ such that
$\left|\frac{\partial g^{il}}{\partial x_j}(x)\right|$ is bounded
and, for $|x|>>0$, there exists $C>0$ such that
\begin{enumerate}
\item for every $i,j$, \begin{equation} \label{1}|g^{ij}(x) - \delta^{ij}
| < \frac{C}{|x|^k}\end{equation}
for some $k>n$;
\item for every $i,j,k,l$
\begin{equation} \label{2} \left|\frac{\partial g_{il}}{\partial x_j}
\right|<
\frac{C}{|x|^k}\end{equation}
\begin{equation} \label{3} \left|\frac{\partial^2 g_{il}}{\partial x_j
\partial
x_k} \right|< \frac{C}{|x|^k}\end{equation} for some $k>n$.
\end{enumerate}
 Then the Laplace-Beltrami operator $\Delta_g$ acting on $1$-forms has
absolutely continuous spectrum equal to $[0,+ \infty)$:
$$[0,+ \infty)= \sigma_{ac}(H_1)=\sigma(H_1).$$
In particular, it has no discrete spectrum. (There might be singularly
continuous spectrum or embedded eigenvalues.)
\end{thm}
\par \noindent
 The main tool for the proof is Belopol'ski-Birman theorem (see
\cite{Reed-Simon}, \cite{Belopol'ski}), which provides a
sufficient condition so that two selfadjoint operators have the
same absolutely continuous spectrum. \par We recall it briefly:

\begin{thm}\label{belopolski} Let $H_0$, $H_1$ be selfadjoint operators
acting respectively on Hilbert spaces ${\mathcal{H}}_0$,
${\mathcal{H}}_1$, and let $E_{\Omega}(H_0)$, $E_{\Omega}(H_1)$, for
$\Omega \subset {\bf R}$, be the associated spectral measures.
\par \noindent If $J\in {\mathcal{L}}({\mathcal{H}}_0,{\mathcal{H}}_1)$
satisfies the conditions:
\begin{enumerate}
\item $J$ has a bounded two-sided inverse;
\item for every bounded interval $I \subset {\bf R}$,
\begin{equation}\label{trace} E_I(H_1)(H_1J-JH_0)E_I(H_0) \in
{\mathcal{I}}_1({\mathcal{H}}_0,{\mathcal{H}}_1),\end{equation}
where ${\mathcal{I}}_1({\mathcal{H}}_0,{\mathcal{H}}_1)= \{ A\in
{\mathcal{L}}({\mathcal{H}}_0,
{\mathcal{H}}_1) \mid (A^*A)^{1/2} \in {\mathcal{I}}_1({\mathcal{H}}_0)\}$ and
${\mathcal{I}}_1
({\mathcal{H}}_0)$ denotes, as usual, the set of trace-class operators on
${\mathcal H}_0$;
\item for every bounded interval $I\subseteq {\bf R}$, $(J^*J -I)E_I(H_0)$
is compact;
\item $JQ(H_0)=Q(H_1),$
where $Q(H_i)$ is the form domain of the operator $H_i$, for $i=0,1$,
\end{enumerate}
then the wave operators $ W^{\pm}(H_1,H_0;J)$ exist, are complete,
and are partial isometries with initial space $P_{ac}(H_0)$ and
final space  $P_{ac}(H_1)$, where $P_{ac}(H_i)$ denotes, as usual,
the absolutely continuous space of $H_i$, for $i=0,1$.\par As a
consequence, the absolutely continuous spectra of $H_0$ and $H_1$
do coincide.
\end{thm}
\begin{rem}\label{remark}{\em We recall (see \cite{Kato}) that if $H$ is a
densely defined, essentially selfadjoint, positive operator on a
Hilbert space ${\mathcal H}$ and ${\bf h}$ is the associated quadratic
form, the form domain $Q(\overline{H})$ of the selfadjoint
operator $\overline{H}$ is the domain of the closure $\tilde{{\bf
h}}$ of the form ${\bf h}$, that is to say: $ Q(\overline{H})$ is
the set of those $u\in {\mathcal H}$ such that there exists a sequence
$\left\{ u_n\right\} \subset D(H)$ converging to $u$ in ${\mathcal H}$
such that $$ {\bf h}[u_n - u_m ] \longrightarrow 0$$ as $n,m
\rightarrow +\infty$.}
\end{rem}
\section{Proof of Theorem \ref{euclidean}}
We will prove that, for a suitable $J: L^2_1({\bf R}^n,e)
\rightarrow L^2_1({\bf R}^n,g)$, $H_1=\overline{\Delta_g}$,
$H_0=\overline{\Delta_e}$ and $J$ satisfy the conditions of
Theorem \ref{belopolski}, for ${\mathcal H}_0=L^2({\bf R}^n,e)$ and
${\mathcal H}_1=L^2({\bf R}^n,g)$.\par \noindent We begin with the
following
\begin{lem}\label{estimates1}
Let $g$ be as in Theorem \ref{euclidean}. Then there exist
$C$,$C_1>0$, $D$,$D_1>0$ such that
\begin{enumerate}
\item for every $x \in {\bf R}^n$
\begin{equation}\label{eq1}C\leq \sqrt{g(x)} \leq C_1 ;\end{equation}
\item for every $x\in {\bf R}^n$, $v$ in the cotangent space at $x$, $T^*_x({\bf R}^n)$
\begin{equation}\label{eq2} D \sum_i v_i^2 \leq g^{ij}(x)v_iv_j \leq D_1
\sum_i v_i^2.\end{equation}
\end{enumerate}
\end{lem}
Proof.  (\ref{eq1}) follows immediately observing that, for every $x$,
 $\sqrt{g(x)}$ is strictly positive and $\sqrt{g(x)} \rightarrow 1$ as
$|x| \rightarrow +\infty$.\par As for (\ref{eq2}), since the
matrix $g^{ij}(x)$, which expresses the Riemannian metric $g$ in
contravariant form, is a continuous function of $x$ and is
positive, its eigenvalues $\lambda_1(x),...,\lambda_n(x)$ depend
continuously on $x$ and are strictly positive. Hence, the
functions $f$ and $h$ defined by $$f(x):= \inf_i \lambda_i (x)$$
and $$h(x):=\sup_i \lambda_i(x),$$ are continuous and strictly
positive. Moreover, since the metric $g$ is asymptotically
euclidean, $f(x)\rightarrow 1$ and $h(x) \rightarrow 1$ as
$|x|\rightarrow +\infty$. As a consequence, there exist $D,D_1>0$
such that, for every $x\in {\bf R}^n$, $$ D \leq f(x) \leq h(x)
\leq D_1,$$ which yields (\ref{eq2}). \qed
\par \bigskip \noindent Lemma \ref{estimates1} implies that there is a
natural identification between $L^2_1( {\bf R}^n,g)$ and $L^2_1(
{\bf R}^n,e)$, and, moreover, (\ref{L2e}) and (\ref{L2G}) are
equivalent norms. As a consequence, the identity map on
$\Lambda^1_c({\bf R}^n)$ extends to a bounded linear operator $$J:
L^2_1({\bf R}^n,e) \longrightarrow L^2_1({\bf R}^n,g),$$ with
bounded two-sided inverse, and condition 1 of Theorem
\ref{belopolski} is satisfied.
\par \bigskip
\noindent

In order to prove (\ref{trace}), we need two Lemmas:

\begin{lem} \label{trace class}
Let $A: \xi \mapsto A_{\xi}$ be a $n \times n$-matrix-valued
function on ${\bf R}^n$, and let ${\mathcal A}$ be the linear operator
$${\mathcal A}:D({\mathcal A})\subset L^2({\bf R}^n) \oplus ... \oplus
L^2({\bf R}^n) \longrightarrow L^2({\bf R}^n) \oplus ... \oplus
L^2({\bf R}^n)$$ of the form $$f(x)A(-i \nabla_x),$$ where $f(x)$
is a function on ${\bf R}^n$ and $A(-i \nabla_x)$ is the operator
$$A(-i \nabla_x)= {\mathcal F} \circ \hat{A}_{\xi} \circ {\mathcal
F}^{-1},$$
 ${\mathcal F}$ being the Fourier transform and $\hat{A}_{\xi}$
 the multiplication operator
 $$ v \longmapsto \hat{A}_{\xi} v$$
 $$ (\hat{A}_{\xi}v)(\xi)= A(\xi)v(\xi).$$
 Let $L^2_{\delta}({\bf
R}^n)$ be the space of functions $h$ such that
$$\|h\|^2_{\delta}=\|(1+|x|^2)^{\frac{\delta}{2}}h(x)\|_{L^2}<
\infty.$$ If, for some $\delta > \frac{n}{2}$, $f(x) \in
L^2_{\delta}({\bf R}^n)$ and, for every pair of indices
$(\alpha,\beta)$, $A_\alpha^\beta(\xi) \in L^2_{\delta}({\bf
R}^n)$, then ${\mathcal A}$ is a trace-class operator.\end{lem} Proof
of Lemma \ref{trace class}. It suffices to show that, for every
fixed $(\alpha, \beta)$, the operator ${\mathcal A}^{\alpha}_{\beta}$
$${\mathcal A}^{\alpha}_{\beta}:D({\mathcal A}^{\alpha}_{\beta}) \subset
L^2({\bf R}^n)\oplus...\oplus L^2({\bf R}^n) \longrightarrow
L^2({\bf R}^n)\oplus ... \oplus L^2({\bf R}^n)$$
$$\omega=(\omega_1,...,\omega_n) \longmapsto
\left(\underbrace{0,...,0,f(x)A^{\beta} _{\alpha}\left(-i \nabla_x
\right) \omega_{\beta}}_{\alpha} ,0,...0\right)$$ is trace-class.
But this latter operator coincides with the composition
$$I_{\alpha} \circ \left(f(x) A^{\beta}_{\alpha} \left(-i \nabla_x
\right)\right)\circ P_{\beta} ,$$
 where
$P_{\beta} $ is the projection $$P_{\beta}
:\,L^2({\bf R}^n) \oplus ...\oplus L^2({\bf
R}^n)\longrightarrow L^2({\bf R}^n)$$ $$\omega=(\omega_1,...,
\omega_{\beta},...,\omega_n) \longmapsto \omega_{\beta} ,$$
$I_{\alpha} $ is the immersion $$I_{\alpha}:\, L^2({\bf
R}^n)\longrightarrow L^2({\bf R}^n) \oplus ... \oplus L^2({\bf
R}^n)$$ $$\omega\longmapsto (\underbrace{0,...,0,\omega}_{\alpha}
,...,0),$$ and $A^{\alpha} _{\beta} (-i\nabla_x)$ is the operator
$$D(A^{\alpha} _{\beta} (-i\nabla_x)) \subset L^2({\bf R}^n)
\longrightarrow L^2({\bf R}^n)$$ $$A^{\alpha} _{\beta}
(-i\nabla_x)= {\mathcal F}\circ(\hat{A_{\xi}})^{\alpha}
_{\beta}\circ{\mathcal F}^{-1},$$ where ${\mathcal F}$ is the Fourier
transform and $(\hat{A_{\xi}})^{\alpha} _{\beta}$ is the
multiplication operator associated to the scalar function
$(A_{\xi})^{\alpha} _{\beta}$.\par The conclusion follows from the
fact that $P_{\beta}$ and
 $I_{\alpha}$ are bounded operators and
  (see \cite{Reed-Simon}, Theorem
XI.21) any operator $$L^2({\bf R}^n) \longrightarrow L^2({\bf
R}^n)$$ of the form $f(x)\,h(-i \nabla_x)$ is trace-class if
$f(x)$ and $h(\xi)$ belong to $L^2_{\delta}({\bf R}^n)$.  \qed
\par \bigskip

\begin{lem}\label{f}
If $f:{\bf R}^n \longrightarrow {\bf R}$ is continuous and  such
that for some $k>n$
\begin{equation}\label{Cxk} |f(x)| <
\frac{C}{|x|^k}\end{equation} when $|x|>>0$, then $f\in
L^2_{\delta} ({\bf R}^n)$ for some $\delta > \frac{n}{2}$.
\end{lem}
Proof. Choosing $\epsilon >0$ such that $k>n + \epsilon$, then  $$
|f(x)|
< \frac{C}{|x|^k}$$
for $|x|>>0$; as a consequence, a straighforward computation in
polar coordinates shows that, for $\delta = \frac{n}{2} +
\epsilon$, $$ \int_{{\bf R}^n} |f(x)|^2 (1+|x|^2)^{\delta} dx <
+\infty .$$ \qed
\par \bigskip
Now, to prove (\ref{trace}), it suffices to see that for every
bounded interval $I \subset {\bf R}$, $$(H_1 - H_0)E_I(H_0)\in
{\mathcal{I}}_1(L^2({\bf R}^n)\oplus ... \oplus L^2({\bf R}^n)).$$ Let
$\Gamma_{ik}^{\alpha}$ be the Christoffel symbols of the
Riemannian connection induced by $g$; then the difference
$H_1-H_0$ is given by
\begin{equation}\label{difference} ((H_1-H_0)\omega)_k = (-g^{ij} +
\delta^{ij})
\frac{\delta^2\omega_k}{\delta x^i \delta x^j} +
g^{ij}\Gamma^{\alpha}_{jk} \frac{\delta \omega_{\alpha}}{\delta x^i} +
g^{ij}\Gamma^{\alpha}_{ij}\frac{\delta \omega_k}{\delta x^{\alpha}}
\end{equation}
$$+g^{ij}\Gamma^{\alpha}_{ik}\frac{\delta \omega_{\alpha}}{ \delta
x^j}+ g^{ij}\frac{\delta \Gamma^{\alpha}_{jk}}{\delta x^i}
\omega_{\alpha} - g^{ij}\Gamma^{\alpha}_{ij}
\Gamma^{\beta}_{\alpha k} \omega_{\beta} - g^{ij}
\Gamma^{\alpha}_{ik} \Gamma^{\beta}_{j \alpha} \omega_{\beta} +
R^i_k \omega_i.$$ A direct computation shows that conditions
(\ref{1}), (\ref{2}), (\ref{3}), and hypothesis 3 in Theorem
\ref{euclidean} imply that $|g^{ij}\Gamma_{jk}^{\alpha}|$,
$|g^{ij}\Gamma_{ij}^{\alpha}|$, $|g^{ij}\Gamma_{ik}^{\alpha}|$,
$\left|g^{ij} \frac{\partial \Gamma_{jk}^{\alpha}}{\partial
x_i}\right|$, $|g^{ij} \Gamma_{ij}^{\alpha} \Gamma_{\alpha
k}^{\beta}|$, $|g^{ij} \Gamma_{ik}^{\alpha} \Gamma_{\alpha
j}^{\beta}|$, $|R^i_k| $ are all bounded from above by $
\frac{C}{|x|^k}$ for some constant $C>0$ and some $k>n$. \par
$(H_1 - H_0) (E_I(H_0))$ is a sum of operators of type $f(x)A(-i
\nabla_x)$, with $f(x) \in L^2_{\delta}({\bf R}^n)$ in view of
Lemma \ref{f}, and $A(\xi)$ smooth and compactly supported.\par

Thus, thanks to Lemma \ref{trace class}, $(H_1 - H_0) (E_I(H_0))$ is
trace-class and condition 2 is fulfilled.
\par \bigskip
\noindent As for condition 3, first of all we observe that the
adjoint of $J$ $$J^*:L^2_1({\bf R}^n,g) \longrightarrow L^2_1({\bf
R}^n,e)$$  satisfies the equation $$\int_{{\bf R}^n} g^{ij}
\omega_i \phi_j \sqrt{g} \,dx=
 \int_{{\bf R}^n} \delta^{ij} \omega_i (J^*\phi)_j dx,$$
and therefore $$(J^*\phi)_k =  \delta_{ik}g^{ij} \phi_j
\sqrt{g}.$$ As a consequence, in local coordinates
$$((J^*J-I)\phi)_k = (\sqrt{g}g^{jk} - \delta^{jk})\phi_j;$$ now,
 $$\left| \sqrt{g}g^{jk} - \delta^{jk} \right| \leq
\left| \sqrt{g}\, \right| \left| g^{jk}- \delta^{jk}\right| +
|(\sqrt{g}-1)|\,\delta^{ij}.$$ By (\ref{1}), there exists $C>0$
such that $$\left| g^{jk}- \delta^{jk}\right| \leq
\frac{C}{|x|^k}$$ for some $k>n$ for $|x|>>0$; moreover, $$
|(\sqrt{g}-1)| = \frac{1}{2} | 1-g|+
o\left(\frac{1}{|x|^k}\right)\leq \frac{K}{|x|^k}$$ for some $K>0$
and some $k>n$ as $|x|\rightarrow + \infty.$\par \noindent Thus,
$\sqrt{g}g^{jk} - \delta^{jk}$ belongs to $L^2_{\delta}({\bf
R}^n)$. Hence $(J^*J-I)E_I(H_0)$ is an operator of type
$f(x)A(-i\nabla_x)$, with $f(x)$ in $L^2_{\delta}({\bf R}^2)$ and
$A(\xi)$ smooth and compactly supported; by Lemma \ref{trace
class}, it is trace-class, and therefore it is compact.\par
\bigskip \noindent As for condition 4, thanks to Remark
\ref{remark}, $Q(H_0)$ and $Q(H_1)$  can be characterized as
follows:
\begin{lem}$Q(H_0)$ is the set of those
$\omega \in L^2_1({\bf R}^n,e)$ for which there exists a sequence
$\{\omega^{(n)}\} \subset \Lambda^1_c({\bf R}^n)$ such that
\begin{equation}\label{omega0}\omega^{(n)} \longrightarrow \omega \quad
\mbox{in} \quad L^2_1({\bf R}^n,e)  \end{equation}
 and
\begin{equation}\label{cauchy0} {\bf h}_0[\omega^{(n)} - \omega^{(m)}]
\longrightarrow 0 \end{equation} as $n,m \rightarrow +\infty$.

Analogously, $Q(H_1)$ is the set of those $\omega \in L^2_1({\bf
R}^n,g)$ such that there exists $\{\psi^{(n)}\}
\subset\Lambda^1_c({\bf R}^n)$ for which
\begin{equation}\label{omega1} \psi^{(n)} \longrightarrow \omega \quad
\mbox{in} \, L^2({\bf R}^n,g) \end{equation}
and
\begin{equation}\label{cauchy1} {\bf h}_1[\psi^{(n)} - \psi^{(m)}]
\longrightarrow 0 \end{equation} as $n,m \rightarrow +\infty$.
\end{lem}
\par \bigskip
\noindent We prove now that \begin{equation}\label{0c1}Q(H_0)
\subseteq Q(H_1) .\end{equation}
 \noindent
For $\omega \in Q(H_0)$, there exists a sequence $\{\omega^{(n)}\}
\subset \Lambda^1_c({\bf R}^n)$ satisfying (\ref{omega0}) and
(\ref{cauchy0}). Due to the equivalence of the norms (\ref{L2e})
and (\ref{L2G}), $$\omega^{(n)} \longrightarrow \omega \:\:
\mbox{in}\:\: L^2_1({\bf R}^n,g); $$ hence, in order to see that
$\omega \in Q(H_1)$ it suffices to prove that $${\bf
h}_1[\omega^{(n)} -\omega^{(m)}] \longrightarrow 0 $$ as $m,n
\rightarrow +\infty$. \par \bigskip To establish this fact, we
consider first the curvature part of ${\bf
h}_1[\omega^{(n)}-\omega^{(m)}]$,
\begin{equation}\label{curvatureomega}\int_{{\bf R}^n} <R(\omega^{(n)} -
\omega^{(m)}),(\omega^{(n)} - \omega^{(m)})>_g \sqrt{g} \,dx.
 \end{equation}
 The following Lemma holds:
\begin{lem}
There exists $C>0$ such that
\begin{equation}\label{curvature} |\int_{{\bf R}^n} <R \omega, \omega>_g
\sqrt{g}\, dx | \leq C \| \omega \|_{L^2_1({\bf
R}^n,e)}^2\end{equation} for every  $\omega \in L^2_1({\bf
R}^n,e)$.
\end{lem}
Proof. Consider for every $x\in {\bf R}^n$ the quadratic form on
$T^*_x({\bf R}^n)$ $$ \omega \longmapsto g^{\alpha \beta}(x)
R^i_{\alpha}(x) \omega_i \omega_{\beta}= R^{i\beta}(x)\omega_i
\omega_\beta.$$ Since  the matrix $R^{i\beta}(x)$ depends
continuously on $x$, its eigenvalues $\lambda_1(x)$,...,
$\lambda_n(x)$ are continuous functions of $x$.  Hence the
function $$f(x):=\sup_i \lambda_i(x),$$ is continuous. Moreover,
since the metric $g$ is asymptotically euclidean, $f(x)\rightarrow
0$ as $|x|\rightarrow +\infty$. As a consequence, there exists
$C>0$ such that $|f(x)|\leq C$ for every $x \in {\bf R}^n$. This
in turn implies $$ \left|R^{i \beta}(x) \omega_i
\omega_{\beta}\right| \leq C \|\omega\|^2_{T^*_x({\bf R}^n)}$$ for
every $x \in {\bf R}^n$ and for every $\omega \in T^*_x({\bf
R}^n)$, which yields (\ref{curvature}). \qed \par
\bigskip
Since $\{\omega^{(n)}\}$ is a Cauchy sequence, the
preceding Lemma implies that
\begin{equation}\label{curvatureomega0}\int_{{\bf R}^n} <R(\omega^{(n)} -
\omega^{(m)}),(\omega^{(n)} - \omega^{(m)})>_g \sqrt{g} \,dx
\longmapsto 0
 \end{equation}
 as $n,m \rightarrow + \infty$.
\par \bigskip
As for the gradient part of ${\bf
h}_1[\omega^{(n)}-\omega^{(m)}]$,
\begin{equation}\label{gradientomega}\int_{{\bf R}^n}|\nabla
(\omega^{(n)} - \omega^{(m)})|^2_g \sqrt{g}\,dx,\end{equation} we
begin by proving
\begin{lem}
There exist $C,D>0$ such that
\begin{equation}\label{CD} C \int_{{\bf R}^n} |\eta|^2_e\, dx \leq
\int_{{\bf
R}^n}
|\eta|^2_g\, \sqrt{g} dx \leq D \int_{{\bf R}^n}
|\eta|^2_e\, dx \end{equation}
for every smooth, compactly supported tensor $\eta=\eta_{ij}$ of rank
2,
where
$$ |\eta|^2_e=\sum_{i,j=1}^n \eta_{ij}^2$$
and
$$ |\eta|_g^2=g^{ij}(x) g^{kl}(x)\eta_{ik}\eta_{jl}.$$
\end{lem}
Proof. Consider, for every $x\in {\bf R}^n$, the quadratic form $$
{\bf R}^{n^2} \times {\bf R}^{n^2} \longrightarrow {\bf R}$$
defined by $$ {\bf \eta}= \eta_{ij} \longmapsto a(x)[{\bf{\eta}}]=
g^{ij}(x) g^{kl}(x)\eta_{ik}\eta_{jl}.$$ Since this quadratic form
is positive and depends continuously on $x$, its eigenvalues
$\lambda_k(x)$, for $k=1,..., n^2$,  are continuous, positive
functions of $x$. Moreover, $$ a(x)[{\bf{\eta}}] \longrightarrow
\sum_{i,j=1}^n \eta_{ij}^2$$ as $|x| \rightarrow + \infty$,
implying that $\lambda_k(x) \rightarrow 1$, for every
$k=1,...,n^2$, as $|x| \rightarrow + \infty$. As a consequence,
there exist $C, D>0$ such that $$ C |{\bf{\eta}}|^2_e \leq
a(x)[{\bf{\eta}}] \leq D |{\bf{\eta}}|^2_e$$ for every $x \in {\bf
R}^n$, ${\bf \eta}\in {\bf R}^{n^2}$, which implies (\ref{CD}).
\qed
\par \bigskip

Setting
 $$ \eta_{ik}= \nabla_i (\omega^{(n)}_k -
\omega^{(m)}_k),$$ (\ref{CD}) yields
\begin{equation} \label{dis1} \int_{{\bf R}^n} |\nabla
(\omega^{(n)} - \omega^{(m)})|_e^2\, dx \leq \frac{1}{C} \int_{{\bf
R}^n} |\nabla (\omega^{(n)} - \omega^{(m)}) |^2_g\, \sqrt{g}\, dx
\end{equation}
and
\begin{equation} \label{dis2} \int_{{\bf R}^n} |\nabla (\omega^{(n)} -
\omega^{(m)}) |^2_g\, \sqrt{g} dx \leq  D \int_{{\bf R}^n}  |
\nabla (\omega^{(n)} - \omega^{(m)})|^2_e \, dx. \end{equation}
Now, $$ \nabla_i \omega_k= \frac{\partial \omega_k}{\partial x_i}
- \Gamma_{ik}^{\alpha}\omega_{\alpha},$$ whence an easy
computation shows that for every $i,k=1,...,n$ $$\left\| \nabla_i
\left( \omega^{(n)}_k - \omega^{(m)}_k \right) \right\|_{L^2({\bf
R}^n,e)} \leq  \left\| \frac{\partial ( \omega^{(n)}_k -
\omega^{(m)}_k)}{\partial x_i}\right\|_{L^2({\bf R}^n,e)} +$$ $$ +
K \left\|\omega^{(n)} - \omega^{(m)}\right\|_{L_1^2({\bf
R}^n,e)}.$$\par

As a consequence, $$ {\bf h}_1[\omega^{(n)} - \omega^{(m)}]
\longrightarrow 0$$ as $n,m \rightarrow 0$. Thus, $Q(H_0)
\subseteq Q(H_1)$.
\par \bigskip
\noindent We complete the proof of Theorem \ref{euclidean} showing
that $Q(H_1) \subseteq Q(H_0)$.\par
\bigskip \noindent
\par \bigskip \noindent For any $\omega \in Q(H_1)$ there exists a
sequence $\{\psi^{(n)}\} \subset \Lambda^1_c({\bf R}^n)$ such that
(\ref{omega1}) and (\ref{cauchy1}) hold.\par Thanks to the
equivalence of the norms (\ref{L2e}) and (\ref{L2G}), $$\psi^{(n)}
\longrightarrow \omega \:\: \mbox{in}\:\: L^2_1({\bf R}^n,e). $$
Thus, in order to see that $\omega \in Q(H_0)$ it suffices to
prove that $$ {\bf h}_0[\psi^{(n)} - \psi^{(m)}] \longrightarrow
0$$ as $m,n \rightarrow +\infty$.\par Now, (\ref{omega1}) and
(\ref{cauchy1}), together with (\ref{curvatureomega}), imply that
$$ \int_{{\bf R}^n} | \nabla ( \psi^{(n)} - \psi^{(m)})|^2_g \,
\sqrt{g}\, dx \longrightarrow 0$$ as $n,m \rightarrow
+\infty$.\par For every $i,k=1,...,n$

 $$ \left\| \frac{\partial (
\psi^{(n)}_k - \psi^{(m)}_k)}{\partial x_i}\right\|_{L^2({\bf
R}^n,e)} \leq$$

$$\leq \left\| \nabla_i \left( \psi^{(n)}_k - \psi^{(m)}_k \right)
\right\|_{L^2({\bf R}^n,e)} + \left\| \Gamma_{ik}^{\alpha}
(\psi^{(n)}_{\alpha} - \psi^{(m)}_{\alpha})\right\|_{L^2({\bf
R}^n,e)} \leq$$

 $$ \leq \left\| \nabla_i
\left( \psi^{(n)}_k - \psi^{(m)}_k \right)\right\|_{L^2({\bf
R}^n,e)} + C \left\|\psi^{(n)} - \psi^{(m)}\right\|_{L_1^2({\bf
R}^n,e)},$$

Then, in view of (\ref{dis1}), $${\bf h}_0[\psi^{(n)}- \psi^{(m)}]
\longrightarrow 0$$ as $n,m \rightarrow +\infty$, and $Q(H_1)
\subseteq Q(H_0)$.
\par \bigskip
Therefore,
$$J(Q(H_0))=Q(H_1).$$
\par \bigskip
\begin{rem}
 Theorem \ref{euclidean} holds, more in general, for
$p$-forms, with $p=1,...n$, with arguments following the same
patterns of the ones developed for $p=1$. Indeed, estimates like
(\ref{eq1}), (\ref{eq2}), (\ref{CD}) hold for $p$-forms, showing
that the identification $J=I: L^2_p({\bf R}^n,g) \rightarrow
L^2_p({\bf R}^n,e)$ is continuous with two-sided bounded inverse.
To establish  the validity of conditions 2.,3. of
Belopol'ski-Birman theorem requires replacing $\Delta_g$ with the
Laplace-Beltrami operator on $p$-forms, given by $$(\Delta_{g(p)}
\omega)_{i_1...i_p}= - \sum_{\alpha,\beta} g^{\alpha \beta}
\nabla_{\alpha} \nabla_{\beta}\, \omega_{i_1...i_p} +
\sum_{j,\alpha} R^{\alpha}_{i_j}
\omega_{i_1...\alpha,\hat{i}_j..i_p}+$$ $$- \sum_{j,
l\not=j,\alpha,\beta} R^{\alpha \; \beta \;}_{\;i_j\;i_l}
\omega_{\alpha i_1...\beta\hat{i}_l...\hat{i}_j...i_p},$$ where
$R^{i\;j\;}_{\;k\;l}$ is the Riemann curvature tensor, which
satisfies the condition $|R^{i\;j\;}_{\;k\;l}|<\frac{C}{|x|^k}$
for $|x|>>0$.\par

The quadratic forms ${\bf h}_0$ and ${\bf h}_1$ have now to be
replaced by $${\bf h}_{0(p)}=\int_{{\bf R}^n}
\sum_{i_1,...,i_p,j=1}^n \left(\frac{\partial
\omega_{i_1...i_p}}{\partial x_j} \right)^2 dx $$ and by
 ${\bf h}_{1(p)}$ expressed by
 $$ {\bf h}_{1(p)}[\omega] =\int_{{\bf R}^n} |\nabla
\omega |^2_g \sqrt{g}\, dx + \int_{{\bf R}^n}
<\tilde R \omega, \omega>_g \sqrt{g}\, dx, $$
where $$|\nabla \omega |^2_g=g^{\alpha \beta}
g^{i_1j_1}...g^{i_pj_p} \nabla_{\alpha}\omega_{i_1...i_p}
\nabla_{\beta} \omega_{j_1...j_p} $$ and $$<\tilde R \omega,
\omega>_g=g^{i_1j_1}...g^{i_pj_p} R_{i_j}^{\alpha}
\omega_{i_1...\alpha... i_p} \omega_{j_1...j_p}+$$ $$+
g^{i_1j_1}...g^{i_pj_p}R^{\alpha\;\beta}_{\;i_j\;i_l}
\omega_{\alpha i_1...\beta ...i_p} \omega_{j_1...j_p}.$$ \par Then
from the fact that a condition similar to (\ref{curvature}) holds
for $<\tilde R \omega, \omega>_g$, the proof follows.

\end{rem}

\end{document}